\documentclass[12pt,twoside]{article}

\usepackage{mathptmx}
\usepackage{enumerate}
\usepackage{longtable}

\usepackage{a4wide}
\usepackage[centertags,intlimits]{amsmath}
\usepackage{amsthm}
\usepackage{amsxtra}
\usepackage{amssymb}
\usepackage{float}
\usepackage{graphicx}

\newcommand{\Real}{\mathbb{R}}

\newcommand{\Natu}{\mathbb{N}}

\newcommand{\eps}{\varepsilon}

\DeclareMathOperator{\diag}{diag}

\newsavebox{\mybox}

\theoremstyle{plain}
\newtheorem{Theo}{Theorem}[section]

\theoremstyle{remark}

\theoremstyle{definition}

\usepackage{dcolumn}
\usepackage{rotating}

\newcommand{\F}{\mathcal{F}}
\newcolumntype{.}{D{.}{.}{-1}}
\newcolumntype{d}[1]{D{.}{.}{#1}}
\newcolumntype{/}{D{/}{/}{-1}}

\begin{document}

\title{Analytical And Numerical Approximation\\of Effective Diffusivities\\in The Cytoplasm of Biological Cells}

\author{Michael Hanke and Marry-Chriz Cabauatan-Villanueva\\Royal Institute of Technology\\
School of Computer Science and Communication}

\date{June 11, 2007}

\maketitle

\begin{abstract}
The simulation of the metabolism in mammalian cells becomes a severe problem if spatial distributions must be taken into account. Especially the cytoplsma has a very complex geometric structure which cannot be handled by standard discretization techniques. In the present paper we propose a homogenization technique for computing effective diffusion constants. This is accomplished by using a two-step strategy. The first step consists of an analytic homogenization from the smallest to an intermediate scale. The homogenization error is estimated by comparing the analytic diffusion constant with a numerical estimate obtained by using real cell geometries. The second step consists of a random homogenization. Since no analytical solution is known to this homogenization problem, a numerical approximation algorithm is proposed. Although rather expensive this algorithm provides a reasonable estimate of the homogenized diffusion constant.
\end{abstract}

\tableofcontents

\section{Introduction}

When mammalian cells are exposed to foreign and potentially harmful compounds a series of events takes place. Following uptake the substance is distributed in different intracellular compartments by diffusion, absorption and desorption. The majority of the compound is either dissolved in the aqueous phase, the cytoplasm, or in the lipophilic phase, the membranes.  Parallel to diffusion and absorption/desorption bioactivation/biotransformation by different soluble and membrane bound enzymes takes place. The purpose of biotransformation is to render the substance suitable for excretion.

\begin{figure}
\begin{center}
\includegraphics[width=0.5\textwidth]{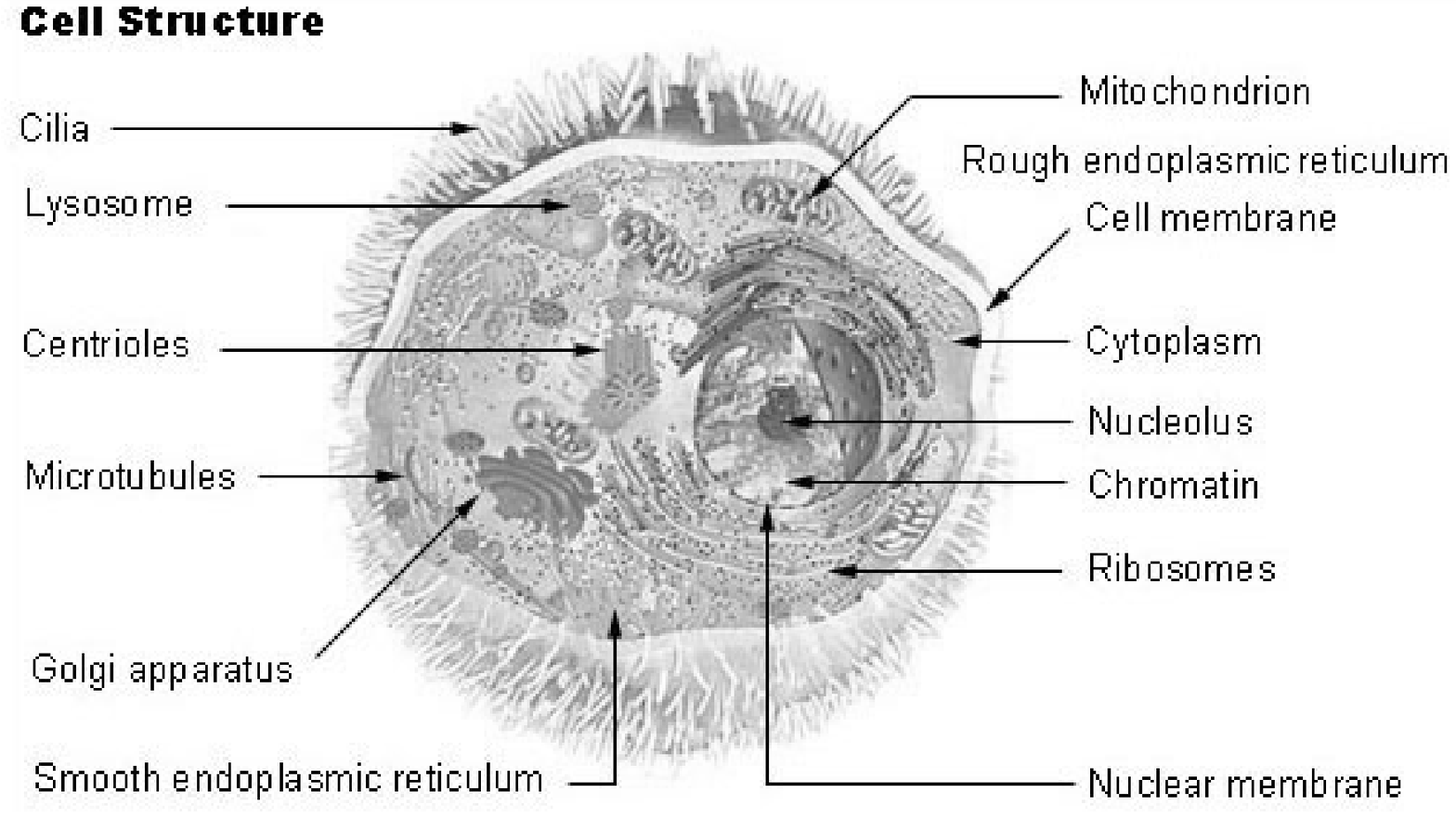}
\end{center}
\caption{Schematic picture of a cell. Picture copyright U.S. National Cancer Institute's Surveillance, Epidemiology and End Results Program, \texttt{http://training.seer.cancer.gov/module\_anatomy/unit2\_1\_cell\_functions\_1.html}\label{fig:schemcell}}
\end{figure}

\begin{figure}
\begin{center}
\includegraphics[width=0.5\textwidth]{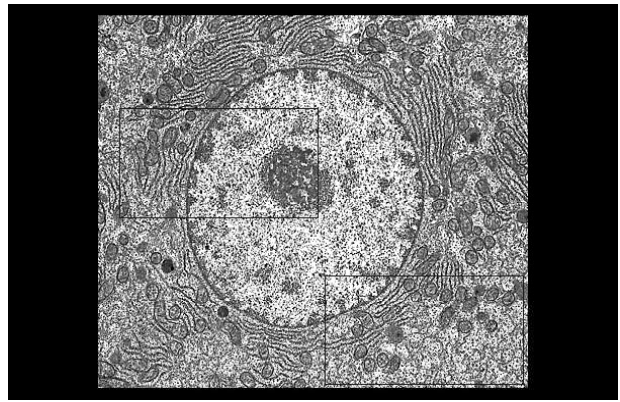}
\end{center}
\caption{Ultrastructure of the cell, nucleus and cytoplasm. Picture copyright Histology Learning System, Boston University, \texttt{http://www.bu.edu/histology/m/index.htm}\label{fig:cellphoto}}
\end{figure}

A human cell consists schematically of an outer cellular membrane, a cytoplasm containing a large number of organelles (mitochondria, endoplasmatic reticulum etc.), a nuclear membrane and finally the cellular nucleus containing DNA. Figure~\ref{fig:schemcell} shows a sketch of a cell while Figure~\ref{fig:cellphoto} shows a microphotograph of a nucleus with part of the surrounding cytoplasm. The organelle membranes create a complex and dense system of membranes or subdomains throughout the cytoplasm. The mathematical description of the biotransformation leads to a system of reaction-diffusion equations in a complex geometrical domain, dominated by thin membranous structures with similar physical and chemical properties. If these structures are treated as separate subdomains, any model becomes computationally very expensive. Moreover, due to the natural variations in the cell structures, every individual cell needs its own mathematical model.

In order to make the system numerically treatable while at the same time retaining the essential features of the metabolism under consideration, in \cite{DrMoJeHa10} a way of homogenizing the cytoplasm has been developed, aiming at a manageable system of reaction-diffusion equations for the various species. In the present paper, we report about numerical experiments which justify some of the strategies in the cited paper. The general modelling assumptions are summarized below. We will use them also in the present report.

Modelling assumptions:
\begin{itemize}
\item On a small scale in space, the volume between the outer cellular membrane and the nucleus membrane consists of layered structures cytoplasm/membranes.
\item In the large scale, this volume contains an unordered set of the small-scale substructures which are uniformly distributed over the volume.
\item The physical and chemical properties of the cytoplasm and of the membranes are uniform.
\item We adopt the continuum hypothesis, i.e., we assume that the set of molecules in the cell can be modelled by considering a continuous representation (a concentration).
\item The processes of absorption and desorption of the individual species into or out of the membrane is much faster than the diffusion and reaction processes. In this case, the relations of the concentrations of a species near a membrane/cytoplasm boundary can be conveniently described with the help of a partition coefficient.
\end{itemize}

In Section~\ref{sec:mathmod}, we introduce the mathematical model. The next section is devoted to a description of the general experimental set-up which is used to compute effective diffusion coefficients numerically. For the solution of the arising boundary-value problems for partial differential equations we used the Comsol Multiphysics\footnote{Comsol Multiphysics is a registered trade mark of Comsol AB, Stockholm, Sweden.} \cite{Co07} environment.

In \cite{DrMoJaHa10}, the diffusion coefficients in the membrane structures have been homogenized by assuming the membranes and the aqueous volumes to be ideal infinite plane layers. This allows for an analytic computation of the effective diffusivity. In Section~\ref{sec:layer}, we will compare effective diffusivities obtained this way against numerically determined effective diffusivities by using computational domains which have been discretized from microphotographs of cell membranes.

The result of the first homogenization step leads to anisotropic diffusion tensors valid locally. Invoking the assumption about the random distribution of the orientation of the membranes, the next step consists of a stochastic homogenization. In contrast to the one- and two-dimensional case, no analytical solution in the general three-dimensional case is known. We will compute the effective diffusion coefficient numerically by Monte Carlo techniques in Section~\ref{sec:rand}.

\section{The Mathematical Model\label{sec:mathmod}}

\subsection{The Governing Equations}

We intend to derive a homogenized model of the reaction and diffusion processes inside the cytoplasm. For that purpose, let $G$ denote the volume between the outer cell membrane and the nucleus membrane (excluding the membranes themselves). This volume is split into two disjoint parts $G_l$ and $G_w$ which denote the lipophilic part and the aqueous part, respectively, of the cell. Note that these subdomains are not necessarily connected. Assume that we are interested in the contrations $c_1,\ldots,c_n$ of $n$ species inside of $G$. For the $k$-th species, it holds
\begin{equation} \label{sreactdiff}
\frac{\partial}{\partial t}c_k=\nabla\cdot(d_k(x)c_k)+R_k(c_1,\ldots,c_n,x),
\quad x\in G,\quad k=1,\ldots,n.
\end{equation}
Here, $d_k$ denotes the diffusion tensor of the $k$-th species which is assumed to be constant in both $G_l$ and $G_w$. $R_k$ denotes the reaction term. It varies strongly with $x$. In the lipophilic part, $R_k\equiv 0$ because no reactions are taking place there. The concentrations of some of the species can be assumed to be constant over time. As a consequence, many of the reaction terms will be linear.

The partition coefficient, $K_{p,k}$ is the equilibrium ratio of the concentrations of species $k$ between the aqueous compartment and its adjacent lipid compartment. This gives rise to boundary conditions,
\begin{equation}
\label{pcc}
c_{k,w}=K_{p,k}c_{k,l}\quad x\in G_w\cap G_l,
\end{equation}
on the inner boundaries
where $c_{k,w}$ and $c_{k,l}$ denote the concentrations in the aqueous and lipid parts, respectively.

The system (\ref{sreactdiff}) with inner boundary conditions (\ref{pcc}) will be supplemented by (outer) boundary conditions and initial conditions. For the purposes of this paper the precise structure of these conditions is not important.

\subsection{The Model Problem}

Let $G\subset\Real^3$ be a bounded domain which will be splitted into two (not necessarily connected) subdomains $G_1$ and $G_2$ such that
\[
G_1\cap G_2=\emptyset,\quad \overline{G}=\overline{G}_1\cup\overline{G}_2.
\]
The interior boundary will be denoted by $\Gamma$,
\[
\Gamma=\overline{G}_1\cap\overline{G}_2.
\]
Consider the equations
\begin{equation} \label{eqeq}
\frac{\partial}{\partial t}v_i
-\nabla\cdot(d_i(x)\nabla v_i)+r_i(x)v_i=f_i(x),\quad x\in G_i,\quad i=1,2.
\end{equation}
Assume additionally boundary conditions on $\partial G$ and initial conditions on $G$ be given.

On the inner boundary $\Gamma$ the flux must be continuous. Let $\mathbf{n}_i$ denote the outer normal at the boundary of $G_i$. The continuity conditions reads now
\begin{equation} \label{gli}
d_1\frac{\partial v_1}{\partial \mathbf{n}_1}
+d_2\frac{\partial v_2}{\partial \mathbf{n}_2}=0,\quad x\in\Gamma.
\end{equation}
The presence of a partion coefficient between the two phases gives rise to the boundary condition
\begin{equation} \label{ibc}
v_1=K_pv_2,\quad x\in\Gamma.
\end{equation}
This problem can be reduced to a problem in a more standard form by introducing
\begin{equation}\label{defu}
u(x):=\begin{cases} v_1(x), & x\in G_1, \\ K_pv_2(x), & x\in G_2.
\end{cases}
\end{equation}
For this new function $u$, the inner boundary conditions become
\[
\left.
\begin{gathered}
d_1\left.\frac{\partial u}{\partial\mathbf{n}_1}\right|_{G_1}+
\frac{1}{K_p}d_2\left.\frac{\partial u}{\partial\mathbf{n}_2}\right|_{G_2}
=0, \\
\left.u\right|_{G_1}=\left.u\right|_{G_2},
\end{gathered}
\right\}\quad x\in\Gamma.
\]
This motivates the definitions
\begin{gather*}
d=\begin{cases} d_1, & x\in G_1, \\ d_2/K_p, & x\in G_2, \end{cases} 
\quad
\sigma=\begin{cases} 1, & x\in G_1, \\ 1/K_p, & x\in G_2,\end{cases} \\
r=\begin{cases} r_1, & x\in G_1, \\ r_2/K_p, & x\in G_2, \end{cases}
\quad
f=\begin{cases} f_1, & x\in G_1, \\ f_2, & x\in G_2. \end{cases}
\end{gather*}
With these definitions, the problem (\ref{gli}), (\ref{ibc}) becomes equivalent to
\begin{equation} \label{gl}
\sigma\frac{\partial}{\partial t}u-\nabla\cdot(d\nabla u)+ru=f,\quad x\in G
\end{equation}
subject to correspondingly modified initial and boundary conditions.

For later use, let
\begin{equation}
p_1=\frac{|G_1|}{|G|},\quad
p_2=1-p_1=\frac{|G_2|}{|G|}.
\end{equation}

\subsection{Going From The Smallest to The Medium Scale: Homogenization of a Periodic Structure\label{perstruc}}

The homogenization procedure for an equation of the type (\ref{gl}) is proved in \cite{PePeSvWy93}. We cite the basic facts. Consider the following problem,
\begin{equation}
\label{pareps}
\begin{gathered}
\sigma^\eps\frac{\partial}{\partial t}u^\eps+A^\eps u^\eps=f^\eps,\quad
u^\eps(0)=u_0 \\
u^\eps\in L^2(0,T;H^1_0(G)).
\end{gathered}
\end{equation}
Here, the operator $A^\eps$ is given by
\[
A^\eps u := -\frac{\partial}{\partial x_i}\left(d^\eps_{ij}
\frac{\partial}{\partial x_j}u\right)+r^\eps u.
\]
For convenience, we use the Einstein summation convention: If an index appears twice in a multiplicative expression, this expression is understood to implicitly represent the sum over this expression where the index varies between 1 and 3 (the dimension of $G$). Moreover, we assume the following construction of the coefficients:
\[
\sigma^\eps(x)=\sigma(x/\eps),\quad r^\eps(x)=r(x/\eps),\quad
d^\eps_{ij}(x)=d_{ij}(x/\eps),\quad i,j=1,2,3.
\]
The functions $\sigma^\eps$ and $r^\eps$ are assumed to belong to $L^\infty(G)$, and
\[
\sigma\geq\sigma_0>0,\quad r(x)\geq 0  \mbox{ a.e. in } G
\]
for some $\sigma_0\in\Real$.
The functions $d^\eps_{ij}$ are assumed to be measurable and to satisfy the conditions $d^\eps_{ij}=d^\eps_{ji}$ and
\[
\alpha|\xi|^2\leq d^\eps_{ij}\xi_i\xi_j\leq\beta|\xi|^2, \mbox{ a.e. in $G$ for all $\xi\in\Real^3$ and $0<\alpha\leq\beta<\infty$.}
\]
Finally, assume $u_0\in L^2(G)$ and $f^\eps\in L^2(0,T;L^2(G))$.

In order to find the homogenized equation, assume
\[
f^\eps\longrightarrow f\mbox{ weakly in $L^2(0,T;L^2(G))$}.
\]
Let $Y$ be an axis parallel hexahedron in $\Real^3$, that is,
\[
Y=\overset{3}{\underset{i=1}\times} (a_i,b_i).
\]
For a $Y$-periodic function $f$, the mean value is given by
\[
\langle f\rangle := \frac{1}{|Y|}\int_Y f(y)dy.
\]
Assume now that $a_{ij}$, $\sigma$, and $r$ are $Y$-periodic. Then it is possible to consider the problem
\begin{equation}
\label{homprob}
\begin{gathered}
\langle\sigma\rangle\frac{\partial}{\partial t}u+A u=f,\quad
u(0)=u_0 \\
u\in L^2(0,T;H^1_0(G)),
\end{gathered}
\end{equation}
where the operator $A$ is given by
\begin{gather*}
Au = -\frac{\partial}{\partial x_i}\left(d_{\text{eff},ij}\frac{\partial}{\partial x_j}u\right)
+\langle r\rangle u\\
d_{\text{eff},ij} = \left\langle d_{ij}-d_{ik}\frac{\partial\varphi_j}{\partial y_k}
\right\rangle,
\end{gather*}
and $\varphi_j$ is the $Y$-periodic solution of the following local elliptic problem:
\begin{equation}
\label{cellprob}
\begin{gathered}
\frac{\partial}{\partial y_i}\left(d_{ik}(y)
\frac{\partial \varphi_j}{\partial y_k}\right)=
\frac{\partial}{\partial y_i}d_{ij}(y), \\
\varphi_j\in W(Y).
\end{gathered}
\end{equation}
Here, $W(Y)= \{\phi\in H^1(Y)| \varphi\mbox{ is $Y$-periodic and } \langle\varphi\rangle=0\}$.

\begin{Theo} \label{t:hompar}
Under the conditions stated above, (\ref{pareps}) and (\ref{homprob}) have unique solutions $u^\eps\in L^2(0,T;H^1_0(G))$ and $u\in L^2(0,T;H^1_0(G))$, respectively, and it holds
\[
u^\eps\longrightarrow u \mbox{ in $L^2(0,T;H^1_0(G))$ weakly as $\eps\rightarrow 0$}.
\]
\end{Theo}

This theorem is proved in \cite[p. 56]{PePeSvWy93}.\footnote{In the reference, the proof is given for a problem without reaction term. But the proof can easily be generalized.}

In our model problem (\ref{gl}) the cell problems can be simplified considerably. We will assume that, in the smallest scale, aqueous and lipid compartments are perfectly layered. It turns out that in this case the computation for the effective diffusion coefficients leads to a transmission problem which can be solved analytically.

Consider the cell problem (\ref{cellprob}). The material consists of perfect layers $\omega^+$ and $\omega^-$ with thicknesses $a^+$ and $a^-$, respectively. The diffusion coefficients in these layers are $d^+$ and $d^-$, respectively. According to our assumptions,
\[
\frac{a^+}{a^-}=\frac{p_1}{p_2}.
\]
Let us choose the following coordinate: $x_2$ is normal to the interface between $\omega^+$ and $\omega^-$ while $x_1$ and $x_3$ span this interface. The diffusion coefficient $d$ is given by
\[
d_{ij}(x) = \begin{cases} 0, & \mbox{ if $i\neq j$,} \\
d_i^+, & \mbox{ if $i=j$ and $x\in\omega^+$,} \\
d_i^-, & \mbox{ if $i=j$ and $x\in\omega^-$.}
\end{cases}
\]
The cell problem is posed on
\[
Y=(0,l_1)\times(-a^-,a^+)\times(0,l_3).
\]
Then the homogenized diffusivities become (see \cite{Per86})
\begin{gather}
d_{\text{eff},11} = (a^+d_1^++a^-d_1^-)/(a^++a^-),\\
d_{\text{eff},22} = (a^++a^-)/\left(a^+/d_2^++a^-/d_2^-\right), \\
d_{\text{eff},33} = (a^+d_3^++a^-d_3^-)/(a^++a^-),\\
d_{ij} = 0 \mbox{ if $i\neq j$}.
\end{gather}
Note that $d_{\text{eff},11}$ and $d_{\text{eff},33}$ are the arithmetic means while $d_{\text{eff},22}$ is the harmonic mean of both diffusivities $a_i^+$ and $a_i^-$.

\subsection{From The Medium to The Large Scale: Stochastic Homogenization}

In global coordinates, we cannot assume that the coordinate system is oriented in the way that we used above. Consider two Cartesian coordinate systems $(x_1,x_2,x_3)$ and $(z_1,z_2,z_3)$. Assume that a given point $x$ has the representation $z=Tx$ with respect to the $z$-coordinates. Note that $T$ is an orthogonal matrix in that case. Denote the matrix of diffusion coefficients with respect to the $x$-coordinates by $Q^\ast$ and that with respect to the $z$-coordinates by $Q$. Then a short calculation yields
\[
Q^\ast = TQT^{-1} = TQT^T.
\]
This is the point to invoke the next critical assumption: We assume that the volume is tightly packed with substructures of the type considered before, namely layered materials. The key assumption is that all orientations are equally probable. This leads to a stochastic description of the diffusion coefficients. We need a homogenization of operators with random coefficients. A theory for that is provided in \cite{JiKoOl94}.

Note that the mean values $\langle\sigma\rangle$ and $\langle r\rangle$ in (\ref{homprob}) are independent of the orientation of the layers. Therefore, it suffices to consider the stationary diffusion problem
\begin{equation} \label{elleps}
A^\eps u^\eps=f,\quad u^\eps\in H^1_0(G)
\end{equation}
with $G\subset \Real^m$ and
\[
A^\eps u = -\frac{\partial}{\partial x_i} \left(d^\eps_{ij} \frac{\partial}{\partial x_j}u\right)
\]
which is the counterpart of (\ref{pareps}). Assume as before that
\[
d^\eps_{ij}(x) = d_{ij}(x/\eps),\quad i,j=1,\ldots,m.
\]
The randomness of the orientation is modelled by assuming that the matrix $A(y)=\big(d_{ij}(y)\big)$ is statistically stationary with respect to the spatial variable $y\in\Real^m$, or equivalently, that $A(y)$ is a typical realization of a stationary random field.

Let $(\Omega,\F,P)$ be a probability space with $\sigma$-algebra $\F$ and probability measure $P$. Let for each $x\in\Real^m$ a random variable $\xi(x)$ over $(\Omega,\F,P)$ be given. The random field $\xi$ is stationary if it can be represented in the form
\[
\xi(x,\omega)=a(T(x)\omega)
\]
where $a(\cdot)$ is a fixed random variable, $T=T(x):\Omega\rightarrow\Omega$ is a measurable transformation which preserves the measure $P$ on $(\Omega,\F)$. Therefore, for the definition of the coefficients $d_{ij}$ in (\ref{elleps}) it is sufficient to consider a matrix $(d_{ij})$ of random variables $d_{ij}:\Omega\rightarrow\Real$. Realizations of coefficients can then be obtained by setting
\[
d_{ij}(x,\omega)=d_{ij}(T(x)\omega).
\]
Assume in the following that $d_{ij}\in L^\infty(\Omega)$ and
\[
\alpha|\xi|^2\leq d_{ij}(\omega)\xi_i\xi_j, \quad \xi\in\Real^m
\]
for almost all $\omega\in\Omega$ with $\alpha>0$ independent of $\xi$ and $\omega$.

A (deterministic) matrix $d(y)$ is said to admit a homogenization if there exists a constant elliptic matrix $d_{\text{eff}}$ such that for any $f\in H^{-1}(G)$ the solutions $u^\eps$ of the Dirichlet problem (\ref{elleps}) it holds
\begin{align*}
u^\eps \longrightarrow u &\mbox{ in $H^1_0(G)$ weakly as $\eps\longrightarrow0$, and} \\
d^\eps\nabla u^\eps \longrightarrow d_{\text{eff}}\nabla u &\mbox{ in $L^2(G)$ weakly as $\eps\longrightarrow0$,}
\end{align*}
where $u$ is the solution of the Dirichlet problem
\[
-\nabla\cdot(d_{\text{eff}}\nabla u)=f,\quad u\in H^1_0(G).
\]
This definition correspondents to the stationary version of Theorem~\ref{t:hompar}. The following theorem holds true \cite[p. 230]{JiKoOl94}:

\begin{Theo} \label{t:homrand}
Assume additionally that the family of mappings $T(x):\Omega\rightarrow\Omega$, $x\in\Real^m$, forms an ergodic $m$-dimensional dynamical system. Then for almost all $\omega\in\Omega$, the matrix with coefficients $d_{ij}(x)=d_{ij}(T(x)\omega)$ admits homogenization, and the homogenized matrix $d_{\textrm{eff}}$ is independent of $\omega$.
\end{Theo}

Unfortunately, an analytical representation of $d_{\text{eff}}$ is only possible in exceptional cases. We are interested in diffusion coefficients having a representation
\[
d(x)= T(x,\omega)QT(x,\omega)^{-1}
\]
where $T(x)\in SO(m)$ is uniformly distributed in $SO(m)$ and $Q$ is a fixed diffusion tensor. In the two-dimensional case, an analytic solution is provided in \cite[p. 235]{JiKoOl94}. For $m=2$, $d_{\text{eff}}$ is simply a scalar equal to the geometric mean of the eigenvalues of $Q$,
\[
d_{\text{eff}} = \sqrt{\det(Q)}.
\]
Here $\det(Q)$ denotes the determinant of $Q$.

There is no analytical solution known for the case $m=3$.

For later use in the experimental estimation of the effective diffusivity, the following observation is important: According to our assumptions on $d$, the estimate
\[
\|(A^\eps)^{-1}\|\leq \alpha^{-1}
\]
holds true such that, for any $f$ and $l$ in $H^{-1}(G)$,
\[
|\langle l, (A^\eps)^{-1}f\rangle | \leq
\|l\|_{H^{-1}(G)}\alpha^{-1}\|f\|_{H^{-1}(G)}
\]
independently of $\omega\in\Omega$. Consequently, for the expectation values it holds
\begin{equation} \label{expconv}
\mathbb{E}\langle l, (A^\eps)^{-1}f\rangle \longrightarrow 
\mathbb{E}\langle l, A^{-1}f\rangle
\end{equation}
by the dominated convergence theorem.

\section{Numerical Determination of Effective Diffusivities\label{numdet}}

Under the assumption that an effective diffusivity for a given problem exists, the corresponding diffusion constants can be determined experimentally. For that, let $D\subset G$ be a subdomain which is in size comparable to $G$ such that the small scale structure is considerable smaller than the size of $D$. Assume that we want to determine the (scalar) diffusion constant for the diffusion process in $x$-direction. In that case it is convenient to use a cylindrical domain
\[
D=(0,L)\times \omega
\]
with $\omega\subset\Real^2$ being some bounded domain. On $D$ consider the stationary diffusion equation
\[
-\nabla\cdot(d(x)\nabla u)=0, \quad x\in D.
\]
The boundary conditions are selected as follows:
\begin{itemize}
 \item On the boundary $\Gamma_0=\{0\}\times\omega$, a fixed Dirichlet condition is given,
\[
u_{\Gamma_0}=c_0.
\]
\item On the boundary $\Gamma_1=\{1\}\times\omega$, a free diffusion into the surrounding medium is assumed,
\[
-n\cdot(d(x)\nabla u)_{\Gamma_1}=M(u_{\Gamma_1}-c_1).
\]
Here, $M$ is the mass transfer coefficient and $c_1$ is the concentration in the bulk solution outside of $D$.
\item All other boundaries $\Gamma_2=\partial G\setminus(\Gamma_0\cup\Gamma_1)$ are isolated,
\[
 -n\cdot(d(x)\nabla u)_{\Gamma_2}=0.
\]
\end{itemize}
If $d(x)$ would be a constant $d_{\text{eff}}$, it would hold
\begin{align*}
d_{\text{eff}} \frac{c_0-u_{\text{out}}}{L} &= N_{\text{average}}, \\
N_{\text{average}} &= \frac{1}{|\Gamma_1|}\int_{\Gamma_1}M(u-c_1)d\Gamma, \\
u_{\text{out}} &= \frac{1}{|\Gamma_1|}\int_{\Gamma_1} ud\Gamma.
\end{align*}
By $|\Gamma_1|$ we denote the Lebesgue measure of $\Gamma_1$.
If $d(x)$ is varying, these equations can be used as an estimation of the effective diffusivity $d_{\text{eff}}$. In case of an anisotropic  effective diffusivity, the above construction leads to an estimate of the effective diffusivity in $x$-direction, i.e., $d_{\text{eff},11}$.

In the one-dimensional case $m=1$, however, an analytic solution is possible. A simple calculation gives
\[
d_{\text{eff}} = \left(\frac{1}{L}\int_0^L d(x)^{-1} dx\right)^{-1},
\]
which amounts to the harmonic mean.

\section{Theoretical And Experimental Diffusivities For Layered Structures\label{sec:layer}}

The homogenization of layered structures in Section~\ref{perstruc} made use of the assumption that we have ideal planes of different materials with different diffusion tensors. In a real biological cell, this assumption is only approximately fulfilled in small subdomains. Besides the effect of not having the parameter $\eps$ close to zero an additional error is introduced this way. The aim of the present section is to obtain some experimental estimates of how large the error will be. We will start with a real photograph of some cell organelles and extract the geometrical structure of the lipophilic and aqueous layers. Then the diffusivity is estimated using the strategy of Section~\ref{numdet}. This diffusion constant will be compared to the theoretical homogenized value.

\subsection{The Experimental Set-up}

\begin{figure}
 \begin{center}
  \includegraphics[width=0.5\textwidth]{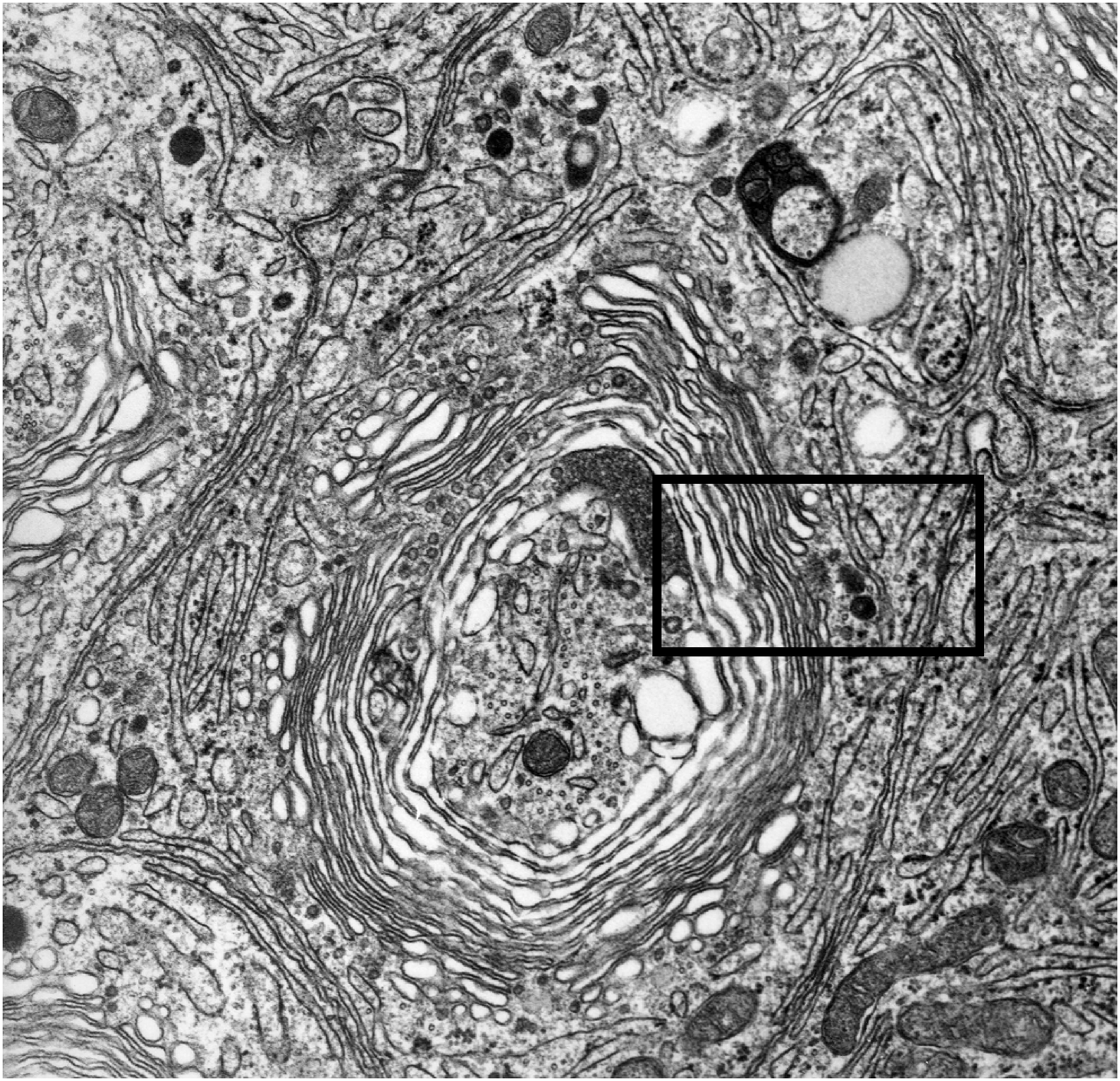}
 \end{center}
\caption{Detail of a rat cell showing the Golgi-apparatus. The box indicates the area used as a reference domain. Copyright Dr.~H.~Jastrow\label{cellphoto}}
\end{figure}

\begin{figure}
 \begin{center}
  \includegraphics[width=0.5\textwidth]{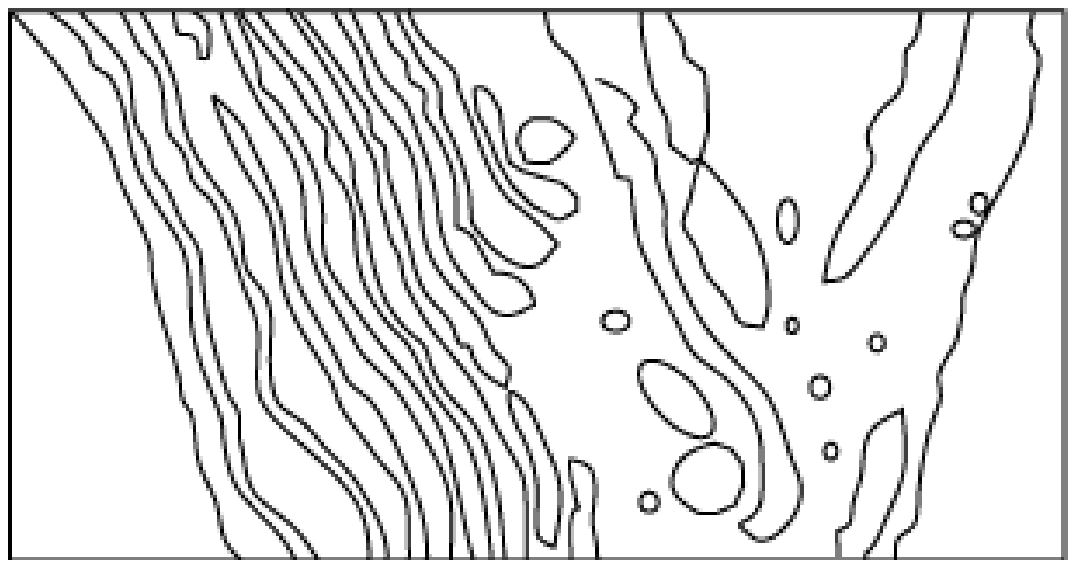}
 \end{center}
\caption{Contrast amplified reference domain. The black areas indicate membranes\label{memb}}
\end{figure}

\begin{figure}
\begin{center}
 \begin{tabular}{c}
   \includegraphics[width=0.6\textwidth]{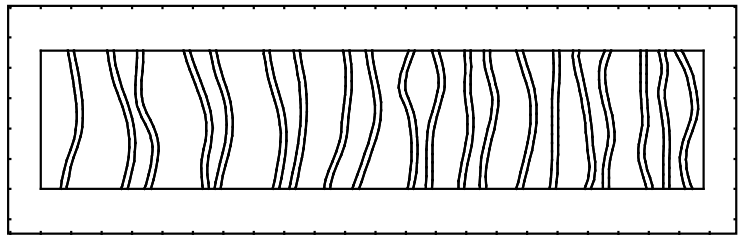} \\
   \textbf{(a)} \\
   \includegraphics[width=1\textwidth]{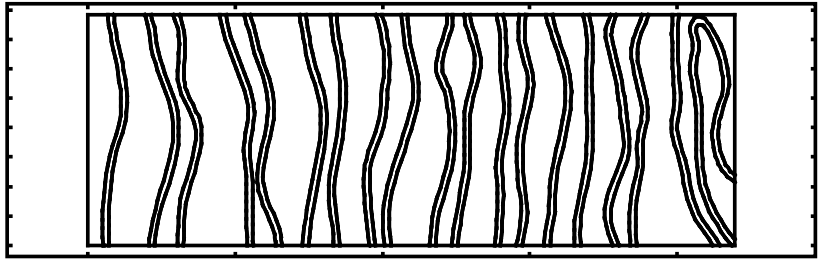} \\
   \textbf{(b)}
 \end{tabular}
\end{center}
\caption{Computational domains for case A \textbf{(a)} and case B \textbf{(b)}\label{compgeo}}
\end{figure}

The experiments of this section are based on the micro-photograph shown in Figure~\ref{cellphoto}. The part enclosed by a box in that figure has been extracted and amplified in contrast. This way, the membrane structure in Figure~\ref{memb} has been obtained. Note that only the black lines represent membranes. The geometry of this structure has been too complex for the software used in the numerical experiments. The number of degrees of freedoms obtained after discretization has become too large. Therefore, we extracted again a part of this geometry in order to make the problem tractable with the available software. Note that the diffusion in this problem is anisotropic. In order to be able to compare the experimental numerical diffusivity with the analytical value, the main orientation of the membranes was aligned with the $y$-axis. The resulting geometries can be found in Figure~\ref{compgeo}. Two cases have been considered.
\begin{itemize}
\item \emph{Case A:} In this case, almost perfect layers have been used.
\item \emph{Case B:} Here we want to estimate the influence of short circuits and more irregular structures.
\end{itemize}

The geometrical data for both data are provided in Table~\ref{geomdata}. The corresponding data for the diffusion constants are given in Table~\ref{diffdata}. Observe that the diffusion in the lipophilic part is anisotropic. This has been used for the numerical experiments. In contrast to that, the homogenized diffusion constant has been determined by using $d_{2,11}$, only. So we expect a larger error in the experiments with the domain of case B.

\begin{table}
\begin{center}
 \begin{tabular}{|l|cc|}
\hline
   & Case A & Case B \\ \hline
$l_x[\textrm{m}]$ & $4.359\times10^{-7}$ & $4.390\times10^{-7}$ \\
$l_y[\textrm{m}]$ & $0.9125\times10^{-7}$ & $1.568\times10^{-7}$ \\
$p_1$ & 0.8122 & 0.8139 \\
$p_2$ & 0.1878 & 0.1861 \\\hline
 \end{tabular}
\end{center}
\caption{Geometric constants\label{geomdata}}
\end{table}

\begin{table}
 \begin{center}
  \begin{tabular}{|l|c|}
\hline
    & Value \\ \hline
$d_1[\textrm{m}^2\textrm{s}^{-1}]$ & $1.0\times 10^{-14}$ \\
$d_{2,11}[\textrm{m}^2\textrm{s}^{-1}]$ & $1.0\times 10^{-12}$ \\
$d_{2,22}[\textrm{m}^2\textrm{s}^{-1}]$ & $1.0\times 10^{-10}$ \\
$K_p$ & $1.26\times 10^{-2}$ \\ \hline
  \end{tabular}
 \end{center}
\caption{Diffusion constants\label{diffdata}}
\end{table}

The experimental determination of the effective diffusivity according to Section~\ref{numdet} can be carried out using two approaches:
\begin{enumerate}
 \item Use the original equation (\ref{eqeq}) subject to the inner boundary conditions (\ref{gli}) and (\ref{ibc}).
\item Use the transformed problem (\ref{gl}) without any inner boundary conditions.
\end{enumerate}
In order to be as close as possible to the original problem we have chosen the first alternative for our experiments. Note that, in the case of $K_p=1$, both approaches are identical.

Unfortunately, it is not possible to formulate the inner transmission conditions (\ref{gli}), (\ref{ibc}) directly in Comsol Multiphysics. Instead, both conditions have been coupled by a penalty approach as suggested in \cite{CoC07}. For a suitably chosen constant $\kappa$, (\ref{gli}), (\ref{ibc}) is replaced by
\begin{equation} \label{penalty}
\begin{aligned}
d_1\frac{\partial v_1}{\partial\mathbf{n}_1} &= \kappa(v_1-K_pv_2),\quad x\in\Gamma \quad (\text{in $G_1$}), \\
d_2\frac{\partial v_2}{\partial\mathbf{n}_2} &= \kappa(K_pv_2-v_1),\quad x\in\Gamma \quad (\text{in $G_2$}).
\end{aligned}
\end{equation}
$\kappa$ acts as a mass transfer coefficient.

For comparison purposes, even the homogenized problem (\ref{homprob}) has been implemented in Comsol Multiphysics.

\subsection{Results}

The experiments have been carried out using the values
\[
 \kappa=10^{-4},\quad M=10^{-7},\quad c_0 = 1, \quad c_1 = 0.
\]
The penalty parameter has been chosen such that both sides of the equations (\ref{penalty}) are somehow in balance. The value of $M$ has been chosen such that the outflow has the order magnitude $0.3c_0$. In case 1, $K_p=1$ while, in case 2, $K_p=0.0126$. The results are summarized in Table~\ref{effdiff}.

\begin{table}
 \begin{center}
  \begin{tabular}{|l|cc.|}
\hline
case & hom. constant & exper. constant & \multicolumn{1}{c|}{rel. difference} \\ \hline
1A & 1.2284 & 1.3131 & 6.9\% \\
1B & 1.2258 & 1.3590 & 10.9\% \\
2A & 1.2312 & 1.2910 & 4.9\% \\
2B & 1.2286 & 1.4680 & 19.5\% \\ \hline
  \end{tabular}
 \end{center}
\caption{Homogenized and experimental effective diffusivities scaled by $10^{-14}$\label{effdiff}}
\end{table}

\begin{figure}
 \begin{center}
  \includegraphics[width=0.8\textwidth]{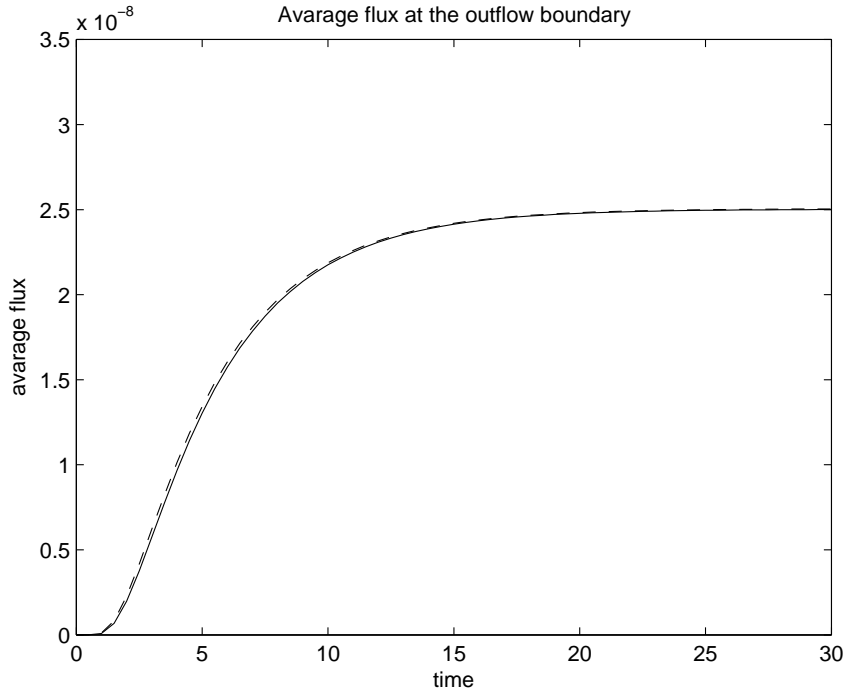}
 \end{center}
\caption{Comparison of the flux at the outflow boundary for the homogenized model (line) and the detailed model (dashed line)\label{transient}}
\end{figure}

The effective diffusivities given above refer to the steady state. In order to get a feeling for the influence of the homogenization on the transient behavior, we compared the time history of the mean flux out of the domain at the left boundary between the original equation (\ref{eqeq}) and its homogenized counterpart (\ref{homprob}). For that, the boundary value problem has been solved as before using the initial condition
\[
 u(x,y) = 10^{-6}\exp\left(-1000\left(\frac{x}{2.179\times10^{-7}}\right)^2\right)\quad \text{ at }t=0.
\]
The value of the experimental effective diffusion has been used in the homogenized problem.
The results for the four different cases differ only marginally. As expected from the experiments, the largest differences occur in case 2B. This is shown in Figure~\ref{transient}.

Summarizing, the following sources of errors occur:
\begin{itemize}
\item The sizes of the sub-structures are not infinitesimal small;
\item The membrane layers are not ideal planes;
\item In the geometry case B, the membranes are touching the outflow boundary;
\item For computing the homogenized diffusion coefficient, only the normal part of the membrane diffusion tensor has been used;
\item The partition coefficients are handled by a penalty approach.
\end{itemize}
The size of the numerical errors is negligible compared to the ones given above.

\section{Experimental Effective Diffusivities in Random Media\label{sec:rand}}

\subsection{The Experimental Set-up\label{expe}}

The idea for estimating the effective diffusivity in the present case is to use a Monte Carlo simulation. For that, the test domain $D$ of Section~\ref{numdet} is chosen to be the unit cube, $D=(0,1)^3$. Let
\[
Q=\diag(d_{11},d_{22},d_{33})
\]
be a fixed diffusion tensor.
 For a given $N\in\Natu$, this cube is subdivided into $N^3$ sub-cubes
\[
D_{ijk} = (x_{i-1},x_i)\times(y_{j-1},y_j)\times (z_{k-1},z_k)
\]
with $x_i=y_i=z_i=ih$ and $h=1/N$. $h$ plays the r\^{o}le of $\eps$ in Theorem~\ref{t:homrand}. One experiment consists of choosing a realization $d^\eps$ such that
\[
\left.d^\eps\right|_{D_{ijk}} = T_{ijk}QT_{ijk}^T
\]
where $T_{ijk}\in SO(3)$ are drawn uniformly distributed in $SO(3)$.

In order to describe the orientation we will use the Euler angles. Any rotation in $SO(3)$ can be described by three angles, the so-called Euler angles. We will use the convention to first rotate around the $x_3$-axis by the angle $\alpha$, then around the (new) $x_1$-axis by $\beta$, and finally around the new $x_3$-axis by $\gamma$. This can be described formally by
\begin{equation}
T=R_3(\gamma)R_1(\beta)R_3(\alpha),\quad \alpha,\gamma\in(0,2\pi),
\quad \beta\in(0,\pi),
\end{equation}
where
\[
R_3(\psi) = \begin{pmatrix} \cos\psi & \sin\psi & 0 \\
                           -\sin\psi & \cos\psi & 0 \\
			       0     &    0     & 1
	    \end{pmatrix},\quad
R_1(\beta) = \begin{pmatrix} 1 &     0     &     0     \\
                             0 & \cos\beta & \sin\beta \\
			     0 &-\sin\beta & \cos\beta
	     \end{pmatrix}.
\]
Let $\mu$ denote the Haar measure on $SO(3)$. Its density has the simple form
\[
d\mu = \frac{1}{8\pi^2}\sin\beta d\alpha d\beta d\gamma
\]
with respect to the Lebesgue measure on $(0,2\pi)\times(0,\pi)\times(0,2\pi)$.

This way, the expectation value of $d_{\text{eff}}^\eps$ can be estimated for given $N$ ($\eps$). The computation of $u_{\text{out}}$ and $N_{\text{average}}$ consists essentially of the evaluation of integrals
\[
\int_{\Gamma_1} ud\Gamma
\]
for functions $u\in H^1(G)$. Since this is a continuous linear functional, we obtain
\[
 d_{\text{eff}}^\eps\longrightarrow d_{\text{eff}} \text{ for }\eps\longrightarrow 0
\]
by using (\ref{expconv}). Since there are no preferred directions in this setting, the diffusion is isotropic such that $d_{\text{eff}}$ is a scalar.

\subsection{Results in 2D}

In the two-dimensional setting, an analytic solution of the random homogenization problem is known. Let
\[
 Q=\diag(d_{11},d_{22}).
\]
The effective diffusivity is the scalar \cite[p. 235]{JiKoOl94}
\[
 d_{\text{eff}} = (d_{11}d_{22})^{1/2}.
\]
We will carry out the experiment described above in the two-dimensional setting in order to obtain a certain gauge for its three-dimensional equivalent.

The two-dimensional counterpart of the experiment described in Section~\ref{expe} is to choose $D=(0,1)^2$ which will be subdivided, for a given $N\in\Natu$, into sub-squares
\[
 D_{ij} = (x_{i-1},x_i)\times (y_{j-1},y_j)
\]
with $x_i=y_i=ih$ and $h=1/N$. The realizations $d^\eps$ are now described by
\[
 \left.d^\eps\right|_{D_{ij}} = T_{ij}QT_{ij}^T
\]
where $T_{ij}\in SO(2)$ are sampled uniformly distributed in $SO(2)$. The elements of $SO(2)$ are simple rotations described uniquely by an angle $\varphi\in[0,2\pi)$,
\[
 T = \begin{pmatrix} \cos\varphi & \sin\varphi \\ -\sin\varphi & \cos\varphi \end{pmatrix}.
\]
The Haar measure $\mu$ on $SO(2)$ has the density $d\mu=\frac{1}{2\pi}d\varphi$ with respect to the Lebesgue measure on $(0,2\pi)$.

The experimental results for
\[
 d_{11}=1,\quad d_{22}=10, \quad d_{\text{eff}}=3.1623
\]
are provided in Table~\ref{randtwod}. We can draw the following conclusions:
\begin{itemize}
 \item The main parameter for the accuracy of the estimation of the effective diffusivity is $N$. This isn't hardly surprising.
\item For a given $N$, the sample size has only a minor influence on the accuracy. Once a certain number of trials has been reached the accuracy does not become better. The optimal sample size seems to be independent of $N$.
\item The standard deviation for sufficiently large sample sizes roughly halves while doubling $N$. This indicates a linear rate of convergence.
\item In all experiments, the mean value of the experimental effective diffusivity is an overestimation of the exact value.
\item If the sample size is too small, the standard deviation is misleading small.
\item The experiments in Section~\ref{sec:layer} indicated an error in the order of magnitude of 5\% -- 10\% between the theoretical homogenized diffusivity and the experimentally observed. These results suggests to use a value of $N=20$ and a sample size of at least 15 trials.
\end{itemize}

\begin{table}
 \begin{center}
  \begin{tabular}{|cc|ccc|}
   \hline
   $N$ & sample size & mean & standard deviation & abs. error \\ \hline
20 & 5 & 3.1693 & 0.0788 & 0.0070 \\
   & 10 & 3.2689 & 0.1604 & 0.1066 \\
   & 15 & 3.1834 & 0.1448 & 0.0211 \\
   & 30 & 3.2225 & 0.1294 & 0.0602 \\
   & 60 & 3.2059 & 0.1569 & 0.0436 \\
   & 90 & 3.1973 & 0.1448 & 0.0350 \\
  & 120 & 3.1708 & 0.1516 & 0.0085 \\
  & 150 & 3.2109 & 0.1371 & 0.0486 \\
  & 180 & 3.1946 & 0.1431 & 0.0323 \\
  & 200 & 3.1971 & 0.1451 & 0.0348 \\ \hline
40 & 5  & 3.2377 & 0.0672 & 0.0754 \\
   & 10 & 3.2272 & 0.0602 & 0.0649 \\
   & 15 & 3.2343 & 0.0675 & 0.0720 \\
   & 30 & 3.2380 & 0.0789 & 0.0757 \\
   & 60 & 3.2496 & 0.0722 & 0.0873 \\
   & 90 & 3.1907 & 0.0746 & 0.0284 \\ \hline
60 & 5  & 3.1950 & 0.0276 & 0.0327 \\
   & 10 & 3.1870 & 0.0487 & 0.0247 \\
   & 15 & 3.1896 & 0.0420 & 0.0273 \\
   & 30 & 3.1916 & 0.0417 & 0.0293 \\ \hline
80 & 15 & 3.2009 & 0.0305 & 0.0386 \\ \hline
  \end{tabular}
 \end{center}
\caption{Experimental effective diffusivities in 2D for $d_{11}=1$, $d_{22}=10$, $d_{\text{eff}}=3.1623$\label{randtwod}}
\end{table}

\subsection{Results in 3D}

Finally, the experimental estimation in the three-dimensional case has been carried out. Unfortunately, the geometry handling in Comsol Multiphysics has led to a severe restriction on how large $N$ can be. Although the machine used had enough memory installed (16 GB), the Java heap space got exhausted rather soon. Moreover, the geometry analysis was surprisingly time-consuming compared to the assembly and solution process.

For the experiments, we have chosen the diffusion coefficients
\[
 d_{11}=9,\quad d_{22}=25,\quad d_{33} =1.
\]
Note that an analytical solution of the homogenization problem is not known. The results are given in Table~\ref{randthreed}.

\begin{table}
 \begin{center}
  \begin{tabular}{|cc|cc|}\hline
$N$ & sample size & mean & standard deviation \\ \hline
 4 &  5 & 7.6753 & 0.4767 \\
   & 10 & 7.2391 & 0.6292 \\
   & 15 & 7.3391 & 0.6667 \\
   & 20 & 7.5785 & 0.8431 \\
   & 30 & 7.5144 & 0.7630 \\ \hline
 8 &  5 & 8.1298 & 0.1226 \\
   & 10 & 8.1251 & 0.2088 \\
   & 15 & 8.0147 & 0.2914 \\
   & 20 & 8.0910 & 0.2395 \\
   & 30 & 8.0490 & 0.2193 \\ \hline
10 &  5 & 8.3499 & 0.1448 \\
   & 10 & 8.2605 & 0.1769 \\
   & 15 & 8.2741 & 0.1523 \\
   & 20 & 8.2783 & 0.1522 \\
   & 30 & 8.3131 & 0.1524 \\ \hline
16 &  5 & 8.6457 & 0.0943 \\
   & 10 & 8.7546 & 0.1005 \\
   & 15 & 8.6834 & 0.0748 \\
   & 20 & 8.6453 & 0.0845 \\
   & 30 & 8.6787 & 0.0752 \\ \hline
20 &  5 & 8.7419 & 0.1162 \\
   & 10 & 8.7383 & 0.0622 \\
   & 15 & 8.7214 & 0.0616 \\
   & 20 & 8.7412 & 0.0505 \\
   & 30 & 8.7281 & 0.0596 \\ \hline
  \end{tabular}
 \end{center}
\caption{Experimental effective diffusivities in 3D for $d_{11}=9$, $d_{22}=25$, $d_{33} =1$\label{randthreed}}
\end{table}

\section{Conclusions}

The present paper explains the homogenization strategy which has been used to derive effective equations for modelling the detailed metabolism in mammalian cells. The cytoplasm has been modelled assuming that three different length scales can be observed. For going from the smallest to the medium scale, an analytic homogenization technique is used. By comparing the analytic effective diffusion constant with results from numerical simulations on real cell geometries taken from photographs an error of 5\% -- 20\% has been observed. Given the accuracy of the known diffusion constants in the lipophilic and aequous parts of the cytoplasm this accuracy appears to be sufficient.

For the step from the medium scale to the large scale, a random homogenization technique has been used. Matheatically the effective diffusivity is known to exists. In the present paper an algorithm has been developed and tested for estimating the homogenized diffusion constant on the large scale. However, the computation times in Comsol Multiphysics have become very large (up to one week on a compute server based on a 2GHz AMD Opteron processor) for a reasonable setup such that alternative solution techniques should be investigated.

The critical assumption in the last step is that about the probability distribution of the structures on the intermediate scale. Its validity can probably only be justified by comparision to biochemical experiments.

More detailed results can be found in \cite{Ca07}.

\paragraph*{Acknowledgement.} The authors are grateful to Dr.~Holger Jastrow, Mainz, for providing us with a high-resolution photograph of cell organelles.

\bibliographystyle{plain}
\bibliography{../../../../art/NEW}

\end{document}